\documentclass[a4paper,11pt]{article}

\usepackage[left=2.5cm,right=2.5cm,top=3cm,bottom=3cm,pdftex]{geometry}
\usepackage{amssymb}
\usepackage{amsmath}
\usepackage{amsthm}

\usepackage{dsfont}
\usepackage{url}
\usepackage{natbib}

\usepackage[utf8]{inputenc}
\usepackage[T1]{fontenc}

\usepackage[pdftex]{graphicx}
\DeclareGraphicsExtensions{.png, .pdf, .jpg} 

\usepackage[pdftex, colorlinks, linkcolor=blue, urlcolor=blue, citecolor=blue, breaklinks=true]{hyperref}

\begin{document}

\begin{center}
\Large\textbf{Optimal Designs for Minimax-Criteria in Random Coefficient Regression Models
} \\[11pt]
\normalsize
Maryna Prus\footnote{Maryna Prus: \href{mailto:maryna.prus@ovgu.de}{maryna.prus@ovgu.de}}\\[11pt]

\footnotesize
Otto-von-Guericke University Magdeburg, Institute for Mathematical Stochastics, 
\\
PF 4120, D-39016 Magdeburg, Germany\\
\normalsize
\end{center}

\begin{quote}
\textbf{Abstract:} We consider minimax-optimal designs for the prediction of individual parameters in random coefficient regression models. We focus on the minimax-criterion, which minimizes the "worst case" for the basic criterion with respect to the covariance matrix of random effects. We discuss particular models: linear and quadratic regression, in detail.

\end{quote}

\section{Introduction}

The subject of this paper is random coefficients regression (RCR) models. These models have been initially defined in biosciences (see e.\,g. \cite{hen3}) and are now popular in many other fields of statistical applications. Besides the estimation of population (fixed) parameters, the prediction of individual random effects in RCR models are often of prior interest. Locally optimal designs for the prediction have been discussed in \cite{pru1} and \cite{pru4}. However, these designs depend on the covariance matrix of random effects. Therefore, some robust criteria like minimax (or maximin), which minimize the largest value of the criterion or maximize the smallest efficiency with respect to the unknown variance parameters, are to be considered. For fixed effects models, such robust design criteria have been well discussed in the literature (see e.\,g. \cite{mue}, \cite{det1}, \cite{sch1}). For optimal designs in nonlinear models see e.\,g. by \cite{paz1}, \cite{pro} and \cite{ell}.

Here we focus on the minimax-criterion for the prediction in RCR models, which minimizes the "worst case" for the basic criterion with respect to the variance parameters. We choose the integrated mean squared error (IMSE) as the basic criterion. We consider particular linear and quadratic regression models in detail.

The structure of this paper is the following: The second part specifies the RCR models and presents the best linear unbiased prediction of the individual random parameters. The third part provides the minimax-optimal designs for the prediction. The paper will be concluded by a short discussion in the last part.

\section{RCR Model}\label{k2}

We consider the RCR models, in which observation $j$ of individual $i$ is given by the following formula:
\begin{equation}\label{mod}
	{Y}_{ij}=\mathbf{f}(x_{j})^\top \mbox{\boldmath{$\beta $}}_i+ \varepsilon_{ij},\quad j=1, \dots, m, \quad i=1, \dots, n, \quad x\in \mathcal{X},
\end{equation}
where $m$ is the number of observations per individual, $n$ is the number of individuals,\linebreak $\mathbf{f} =(f_1, \dots,f_p)^\top$ is a vector of known regression functions. The experimental settings $x_j$ come from an experimental region $\mathcal{X}$.
The observational errors $\varepsilon_{ij}$ are assumed to have zero mean and common variance $\sigma^2>0$. 
The individual parameters $\mbox{\boldmath{$\beta $}}_i=( \beta_{i1}, \dots, \beta_{ip})^\top$ have unknown expected value (population mean) $\mathrm{E}\,(\mbox{\boldmath{$\beta $}}_i)={\mbox{\boldmath{$\beta $}}}$ and known positive definite covariance matrix $\mathrm{Cov}\,(\mbox{\boldmath{$\beta $}}_i)=\sigma^2\mathbf{D}$. 
All individual parameters $\mbox{\boldmath{$\beta $}}_{i}$ and all observational errors $\varepsilon_{ij}$ are assumed to be uncorrelated.

The best linear unbiased predictor for the individual parameter $\mbox{\boldmath{$\beta $}}_i$ is given by
\begin{equation*}
\hat{\mbox{\boldmath{$\beta $}}}_i=(\mathbf{F}^\top \mathbf{F}+\mathbf{D}^{-1})^{-1}(\mathbf{F}^\top \mathbf{F}\,\hat{\mbox{\boldmath{$\beta $}}}_{i;{\rm ind}}+\mathbf{D}^{-1}\hat{\mbox{\boldmath{$\beta $}}}),
\end{equation*}
which is a weighted average of the individualized estimator $\hat{\mbox{\boldmath{$\beta $}}}_{i;{\mathrm{ind}}}=(\mathbf{F}^\top \mathbf{F})^{-1}\mathbf{F}^\top \mathbf{Y}_i$ based only on observations at individual $i$ and the best linear unbiased estimator $\hat{\mbox{\boldmath{$\beta $}}}=(\mathbf{F}^\top \mathbf{F})^{-1}\mathbf{F}^\top \bar{\mathbf{Y}} $ for the population mean parameter.  $\mathbf{Y}_i=(Y_{i1}, \dots,Y_{im})^\top$ is the individual vector of observations, $\bar{\mathbf{Y}} =\frac{1}{n}\sum_{i=1}^n{\mathbf{Y}_i}$ is the mean observational vector and $\mathbf{F}=(\mathbf{f}(x_1), \dots, \mathbf{f}(x_m))^\top$ is the design matrix, which is assumed to be of full column rank.

The mean squared error matrix of the of the vector $\hat{\mathbf{B}}=(\hat{\mbox{\boldmath{$\beta $}}}_1^\top, \dots,\hat{\mbox{\boldmath{$\beta $}}}_n^\top)^\top$ of all predictors of all individual parameters is given by the following formula (see e.\,g. \cite{pru1}):
\begin{equation}\label{mse}
\mathrm{MSE}= \sigma^2\left(\frac{1}{n}\left(\mathds{1}_n\mathds{1}_n^\top \right)\otimes  \left(\mathbf{F}^\top \mathbf{F}\right)^{-1} + \left(\mathds{I}_n-{\textstyle{\frac{1}{n}}}\mathds{1}_{n}\mathds{1}_{n}^\top \right)\otimes \left(\mathbf{F}^\top \mathbf{F}+\mathbf{D}^{-1}\right)^{-1} \right),
\end{equation}
where $\mathds{I}_n$ denotes the identity matrix, $\mathds{1}_n$ is the vector of length $n$ with all entries equal to $1$ and $\otimes$ denotes the Kronecker product.

\section{Optimal Designs}\label{k3}

For this paper we define the exact designs as follows:
\begin{equation*}
\xi= \left( \begin{array}{ccc}  x_1 & , \dots, & x_k \\  m_1 &, \dots,& m_k \end{array} \right),
\end{equation*} 
where $x_1, \dots, x_k$ are the distinct experimental settings (support points), $k\leq m$, and $m_1, \dots, m_k$ are the corresponding numbers of replications. For analytical purposes we will focus on the approximate designs, which we define as
\begin{equation*}
\xi= \left( \begin{array}{ccc}  x_1 & , \dots, & x_k \\  w_1 &, \dots,& w_k \end{array} \right),
\end{equation*}
where $w_j=m_j/m$ and only the conditions $w_j\geq 0$ and $\sum_{j=1}^{k}w_j=1$ have to be satisfied (integer numbers of replications are not required). Further we will use the notation
\begin{equation*}
\mathbf{M}(\xi)=\frac{1}{m}\sum_{j=1}^k m_j\mathbf{f}(x_j)\mathbf{f}(x_j)^\top
\end{equation*} 
for the standardized information matrix from the fixed effects model and $\mbox{\boldmath{$\Delta $}}=m\, \mathbf{D}$ for the adjusted dispersion matrix of the random effects. We assume the matrix $\mathbf{M}(\xi)$ to be non-singular. With this notation the definition of mean squared error matrix \eqref{mse} can be extended for approximate designs  to
\begin{equation*}
\mathrm{MSE}(\xi)= {\frac{1}{n}}\left(\mathds{1}_n\mathds{1}_n^\top \right)\otimes  \mathbf{M}(\xi)^{-1} + \left(\mathds{I}_n-{\frac{1}{n}}\mathds{1}_{n}\mathds{1}_{n}^\top \right)\otimes \left(\mathbf{M}(\xi)+\mathbf{\Delta}^{-1}\right)^{-1}, 
\end{equation*} 
when we neglect the constant term $\frac{\sigma^2}{m}$. 

\subsection{IMSE-criterion}\label{k31}

In this work we focus on the integrated mean squared error (IMSE-) criterion. For the prediction of individual parameters we define the IMSE-criterion (see also \cite{pru1}) as the sum over all individuals
\begin{equation*}
\mathrm{IMSE}_{pred} =\sum_{i=1}^{n}\mathrm{E}\,\left(\int_{\mathcal{X}} (\hat{\mu}_i(x) - \mu_i (x))^2 \nu (\mathrm{d}x)\right)
\end{equation*}
of the expected integrated squared distances of the predicted and the real response, $\hat\mu_i=\mathbf{f}^\top \hat{\mbox{\boldmath{$\beta $}}}_i$ and $\mu_i=\mathbf{f}^\top \mbox{\boldmath{$\beta $}}_i$, with respect to a suitable measure $\nu$ on the experimental region $\mathcal{X}$, which is typically chosen to be uniform on $\mathcal{X}$ with $\nu(\mathcal{X})=1$.
For an approximate design $\xi$ the IMSE-criterion has the form
\begin{equation}\label{phi}
\mathrm{IMSE}_{pred}(\xi) = \mathrm{tr}\left(\mathbf{M}(\xi)^{-1}\mathbf{V}\right)+(n-1)\mathrm{tr}\left(\left(\mathbf{M}(\xi)+\mathbf{\Delta}^{-1}\right)^{-1}\mathbf{V}\right),
\end{equation}
where $\mathbf{V}=\int_{\mathcal{X}} \mathbf{f}(x)\mathbf{f}(x)^\top\nu (\mathrm{d}x)$, which may be recognized as the information matrix for the weight distribution $\nu$ in the fixed effects model.

\subsection{Minimax-criteria}

In this section we consider optimal designs for the prediction in particular RCR models: straight line and quadratic regression. We define the minimax-criterion as the worst case of the IMSE-criterion with respect to the unknown variance parameters. 

We additionally assume the diagonal covariance structure of random effects. Then IMSE-criterion \eqref{phi} will increase with increasing values of variance parameters.
However, if all these parameters will be large, the criterion function will tend to the IMSE-criterion in the fixed effects model (multiplied by the number of individuals $n$). Therefore, we fix some of the variances and consider the behavior of minimax-optimal designs in the resulting particular cases. 

Note that for special RCR, where only the intercept is random, optimal designs for fixed effects models retain their optimality (see \cite{pru1}).

\subsubsection*{Straight line regression}

We consider the linear regression model
\begin{equation}\label{lr}
	Y_{ij}= \beta_{i1}+\beta_{i2}x_j+\varepsilon_{ij} 
\end{equation}
on the experimental regions $\mathcal{X}=[0,1]$ with the diagonal covariance structure of random effects: $\mathbf{D}=\textrm{diag} (d_1, d_2)$, and a small intercept variance: $d_1\rightarrow 0$. For the IMSE-criterion we choose the uniform weighting $\nu=\lambda_{[0,1]}$, which leads to $\mathbf{V}=\int_0^1\mathbf{f}(x)\mathbf{f}(x)^\top\mathrm{d}x$. As proved in \cite{pru3}, ch.~5, IMSE-optimal designs for the prediction in model \eqref{lr} are of the form 
\begin{equation*}
\xi_{w}= \left( \begin{array}{cc}  0 & 1 \\ 1-w & w \end{array} \right),
\end{equation*}
where $w$ denotes the optimal weight of observations at the support point $x=1$. Then we obtain the following form of IMSE-criterion \eqref{phi}:
\begin{equation*}
\mathrm{IMSE}_{pred}(\xi) = \frac{1}{3}\left(\frac{1}{mw(1-w)}+(n-1)\frac{d_2}{1+mwd_2}\right).
\end{equation*}
It is easy to see that this criterion increases with increasing values of the slope variance. The latter property allows us to define the minimax-criterion as follows:
\begin{equation*}
\textrm{IMSE}_{max}(\xi) := \textrm{lim}_{d_2\rightarrow \infty}\textrm{IMSE}(\xi),
\end{equation*}
which results in
\begin{equation*}
\textrm{IMSE}_{max}(\xi) = \frac{1}{3m}\left(\frac{1}{w(1-w)}+(n-1)\frac{1}{w}\right)
\end{equation*}
and leads to the following optimal weight:
\begin{equation*}
w^*_{max} = \frac{n-\sqrt{n}}{n-1}.
\end{equation*}
Figure~1 illustrates the behavior of the optimal design with respect to the number of individuals $n$ for all integer values in the interval $[2,500]$. As we can see in Figure~1, the optimal weight increases with increasing number of individuals. Figure~2 presents the efficiency of the minimax-optimal design $w^*_{max}$ with respect to the locally optimal designs in dependence of the rescaled slope variance $\rho={d_2}/{(1+d_2)}$ for fixed numbers of individuals $n=10$, $n=50$ and $n=500$. For all numbers of individuals the efficiency is high and increases with increasing slope variance.
\begin{figure}[ht]
    \begin{minipage}[]{8.2 cm}
       \centering
       \includegraphics[width=78mm]{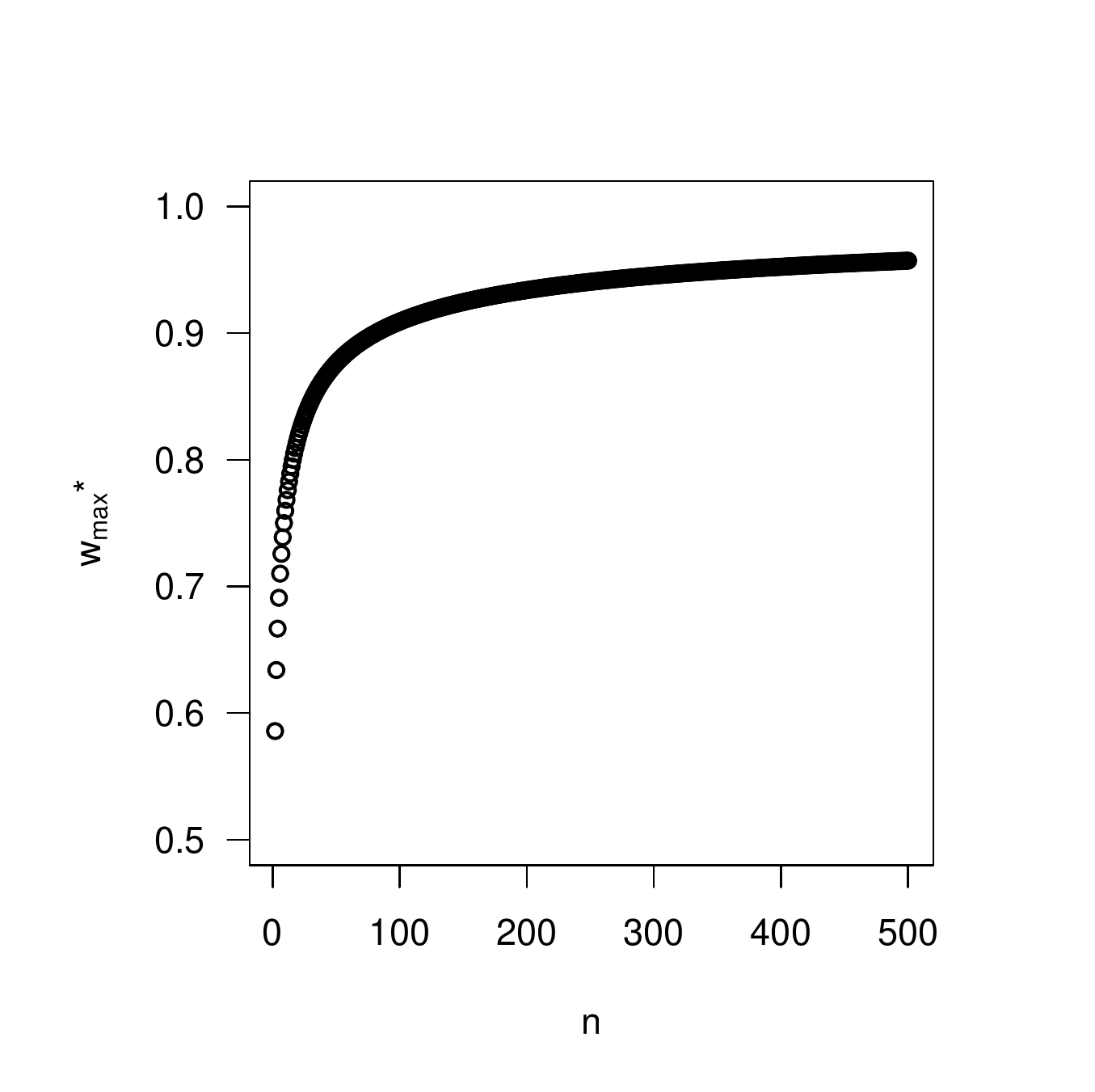}
       \end{minipage}
       \begin{minipage}[]{8.2 cm}
       \centering
       \includegraphics[width=78mm]{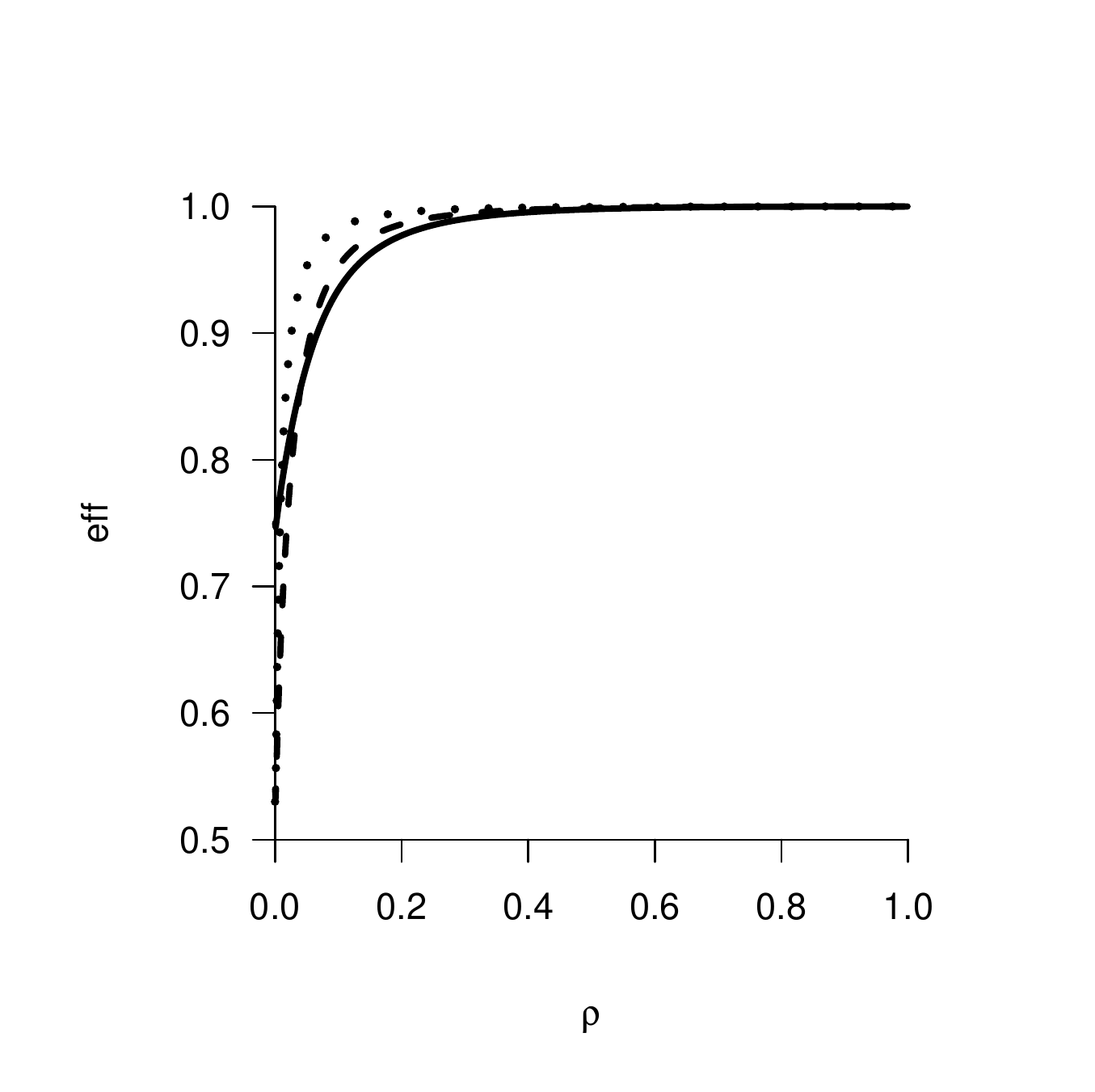}
       \end{minipage}
       \vspace{-1mm}
       \hspace*{5 mm}
       \begin{minipage}[]{6.5 cm}
       Figure 1: Minimax-optimal weight $w^*_{max}$ in dependence of number of individuals $n$ for linear regression
       \end{minipage}
       \hspace{15 mm}
       \begin{minipage}[]{6.5 cm}
       Figure 2: Efficiency of minimax-optimal designs for linear regression
			for $n=10$ (solid line), $n=50$ (dashed line), $n=500$ (dotted line) 
       \end{minipage}
    \end{figure}

\subsubsection*{Quadratic regression}

We investigate the quadratic regression model
\begin{equation}\label{qr}
	Y_{ij}= \beta_{i1}+\beta_{i2}x_j+\beta_{i3}x_j^2+\varepsilon_{ij}
\end{equation}
on the standard symmetric design region $\mathcal{X}=[-1,1]$ with a diagonal covariance matrix of random effects: $\mathbf{D}=\textrm{diag} (d_1, d_2, d_3)$. For the IMSE-criterion we apply the uniform weighting $\nu=\frac{1}{2}\lambda_{[-1,1]}$.

In \cite{pru3}, ch.~5, it has been established that optimal designs in model \eqref{qr}
are of the form
\begin{equation*}
\xi_{w}= \left( \begin{array}{ccc}  -1 & 0 & 1 \\ w & 1-2w & w \end{array} \right).
\end{equation*}

Because of its complexity, the general form of the IMSE-criterion for the quadratic regression will not be presented here. As it was mentioned at the beginning of this section, the IMSE-criterion increases with increasing variances. Hence, we will fix some of the variances by small values and consider minimax-criteria for the resulting particular cases.

\textbf{Case 1.} $d_1\rightarrow 0$ and $d_2\rightarrow 0$

If both the intercept and the slope variances are small, the worst case of the IMSE-criterion is given by its limiting value ($d_3\rightarrow \infty$). We define, therefore, the minimax-criterion in this case as
\begin{equation*}
\textrm{IMSE}_{max}(\xi) := \textrm{lim}_{d_3\rightarrow \infty}\textrm{IMSE}(\xi),
\end{equation*}
which leads to the optimal design
\begin{equation*}
w^*_{max} = \frac{3n+5-2\sqrt{6n+10}}{6(n-1)}.
\end{equation*}
Figure~3 illustrates the behavior of the optimal weight with respect to the number of individuals $n$.

\textbf{Case 2.} $d_1\rightarrow 0$ and $d_3\rightarrow 0$  

If both variances of the intercept and of the coefficient of the quadratic term are small, the minimax-criterion can be defined as
\begin{equation*}
\textrm{IMSE}_{max}(\xi) := \textrm{lim}_{d_2\rightarrow \infty}\textrm{IMSE}(\xi).
\end{equation*}

\textbf{Case 3.} $d_3\rightarrow 0$

If only the variance of the coefficient of the quadratic term is small, we obtain the following worst case of the IMSE-criterion:
\begin{equation*}
\textrm{IMSE}_{max}(\xi) := \textrm{lim}_{d_1, d_2\rightarrow \infty}\textrm{IMSE}(\xi).
\end{equation*}
For both cases 2 and 3 we receive the following optimal weight of the observations at the support point $1$:
\begin{equation*}
w^*_{max} = \frac{5n+3-2\sqrt{10n+6}}{10(n-1)},
\end{equation*}
which is described by Figure~4.
\begin{figure}[ht]
    \begin{minipage}[]{8.2 cm}
       \centering
       \includegraphics[width=78mm]{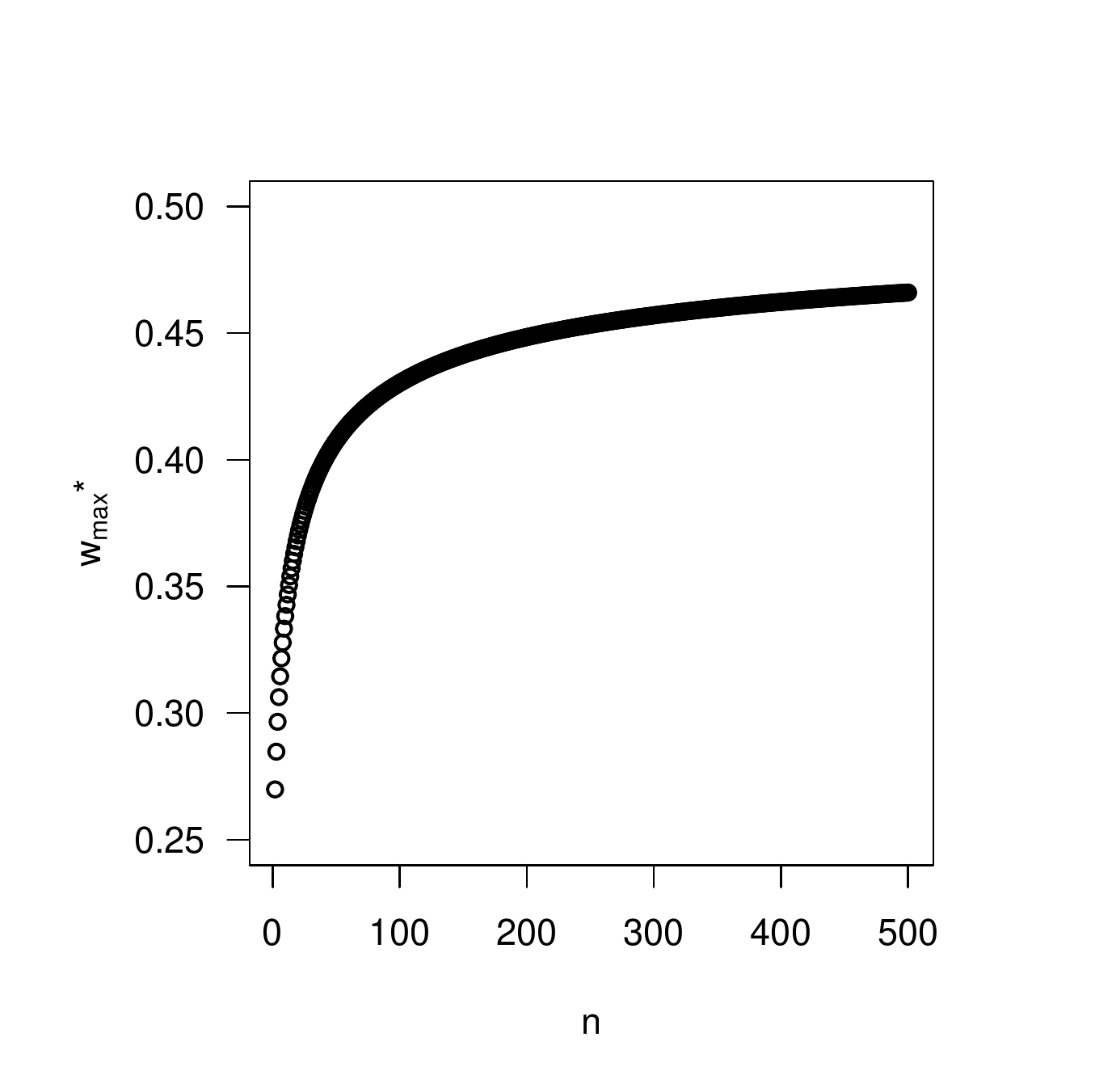}
       \end{minipage}
       \begin{minipage}[]{8.2 cm}
       \centering
       \includegraphics[width=78mm]{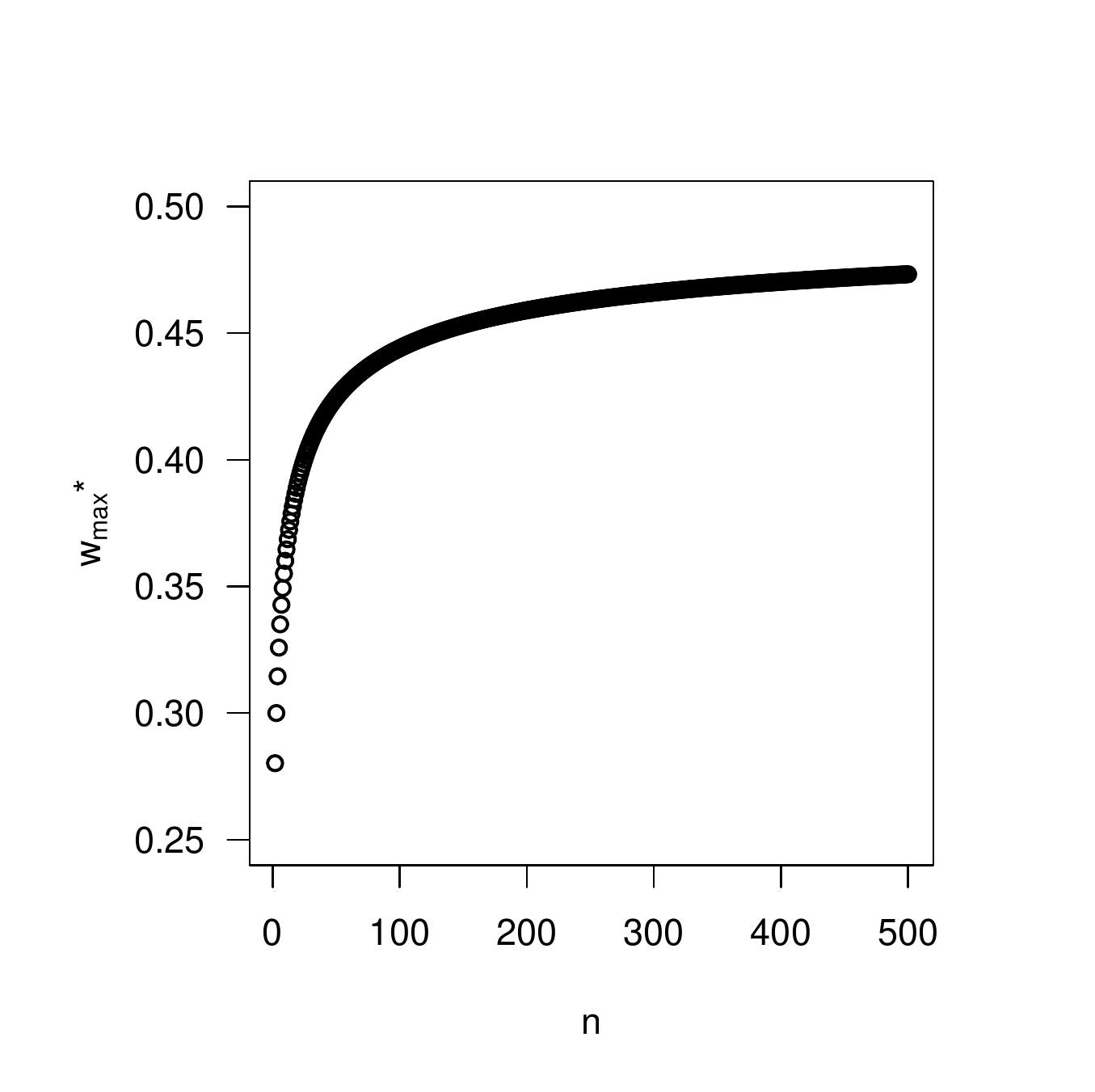}
       \end{minipage}
       \vspace{-1mm}
       \hspace*{5 mm}
       \begin{minipage}[]{6.5 cm}
       Figure 3: Minimax-optimal weight $w^*_{max}$ in dependence of number of individuals $n$ for quadratic regression, case~1
       \end{minipage}
       \hspace{15 mm}
       \begin{minipage}[]{6.5 cm}
       Figure 4: Minimax-optimal weight $w^*_{max}$ in dependence of number of individuals $n$ for quadratic regression, cases~2 and 3
       \end{minipage}
    \end{figure}

\textbf{Case 4.} $d_2\rightarrow 0$

If only the slope variance is small, we determine the minimax-criterion as the limiting value
\begin{equation*}
\textrm{IMSE}_{max}(\xi) := \textrm{lim}_{d_1, d_3\rightarrow \infty}\textrm{IMSE}(\xi)
\end{equation*}
of the IMSE-criterion. The resulting optimal weight is given by the following formula:
\begin{equation*}
w^*_{max} = \frac{-3n-5+2\sqrt{6n^2+10n}}{10(n-1)}.
\end{equation*}

\textbf{Case 5.} $d_1\rightarrow 0$

For small intercept variance the minimax-criterion can be defined as
\begin{equation*}
\textrm{IMSE}_{max}(\xi) := \textrm{lim}_{d_2, d_3\rightarrow \infty}\textrm{IMSE}(\xi),
\end{equation*}
which leads to the minimax-optimal weight
\begin{equation*}
w^*_{max} = \frac{n-\sqrt{n}}{2(n-1)}.
\end{equation*}
The behaviors of the optimal designs in cases~4 and 5 are illustrated by Figures~5 and 6, respectively. 
		\begin{figure}[ht]
    \begin{minipage}[]{8.2 cm}
       \centering
       \includegraphics[width=78mm]{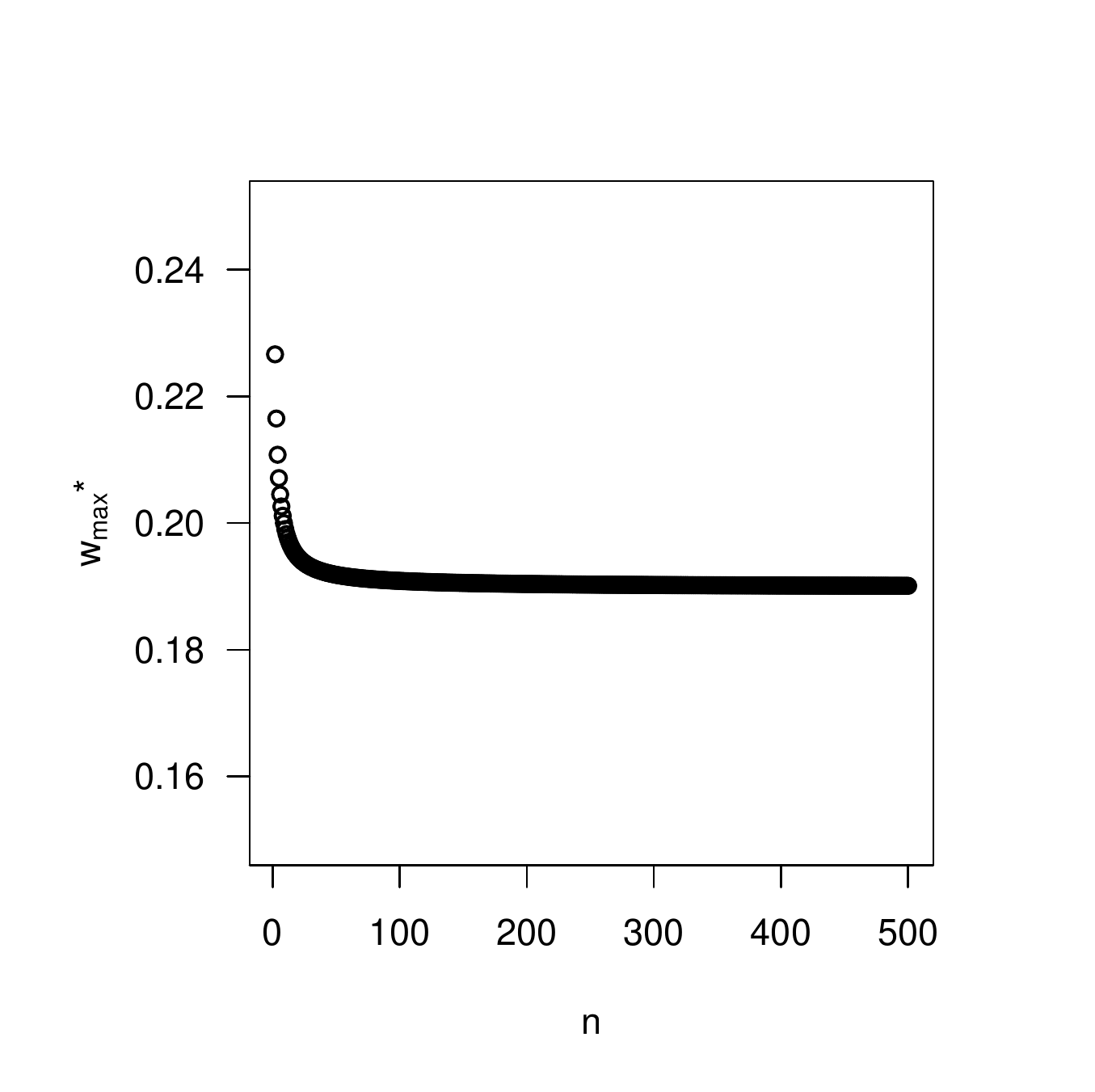}
       \end{minipage}
       \begin{minipage}[]{8.2 cm}
       \centering
       \includegraphics[width=78mm]{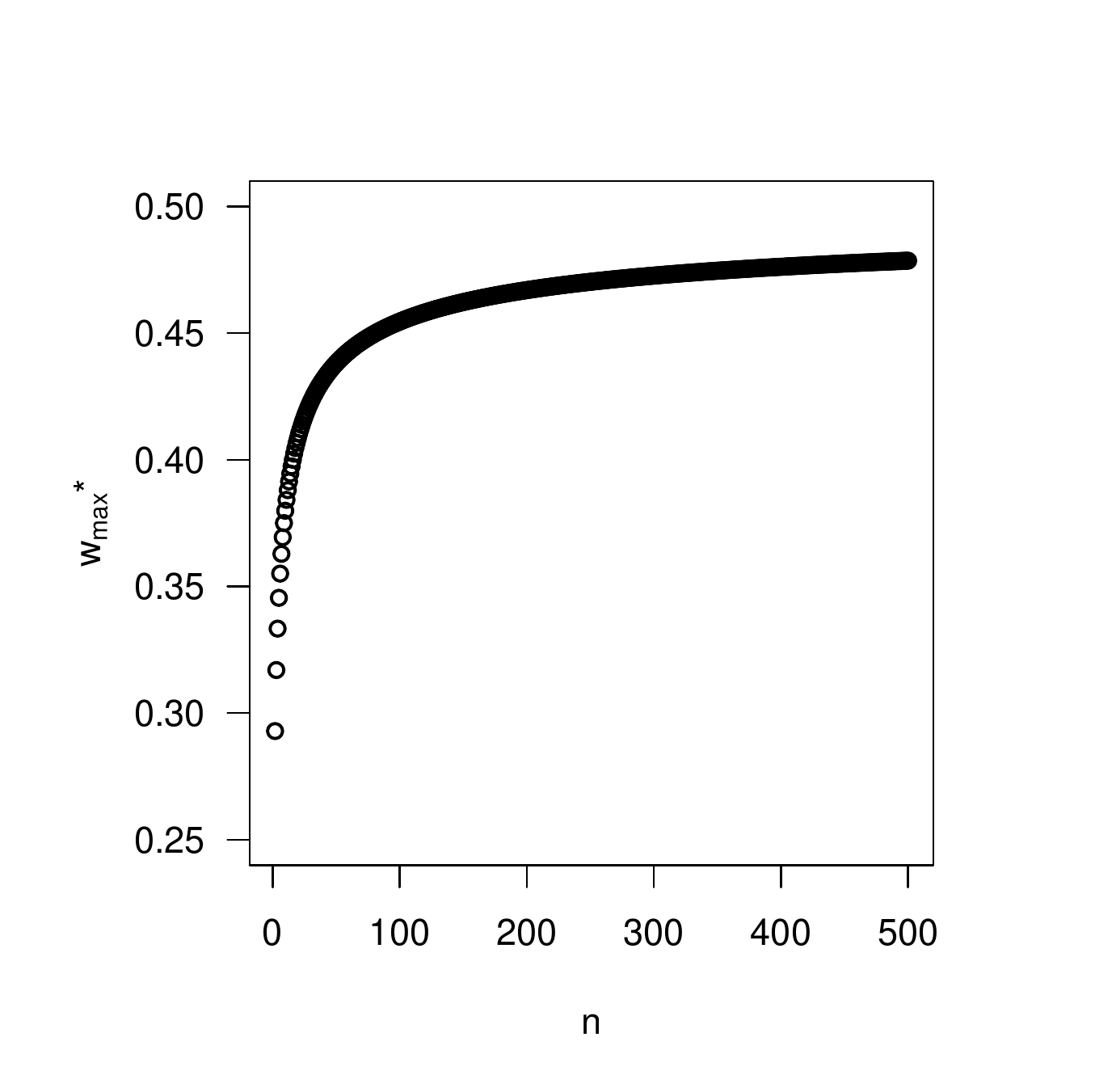}
       \end{minipage}
       \vspace{-1mm}
       \hspace*{5 mm}
       \begin{minipage}[]{6.5 cm}
       Figure 5: Minimax-optimal weight $w^*_{max}$ in dependence of number of individuals $n$ for quadratic regression, case~4
       \end{minipage}
       \hspace{15 mm}
       \begin{minipage}[]{6.5 cm}
       Figure 6: Minimax-optimal weight $w^*_{max}$ in dependence of number of individuals $n$ for quadratic regression, case~5
       \end{minipage}
    \end{figure}
		
As we can see on the graphics, the optimal weights increase with increasing number of individuals $n$ in cases~1, 2, 3 and 5 and decrease in case~4. For cases~1 and 2 we consider  the efficiency of the minimax-optimal designs with respect to the locally optimal designs in dependence of the rescaled variances $\rho={d_3}/{(1+d_3)}$ and $\rho={d_2}/{(1+d_2)}$, respectively, for fixed numbers of individuals (Figures~7 and 8).
		\begin{figure}[ht]
    \begin{minipage}[]{8.2 cm}
       \centering
       \includegraphics[width=78mm]{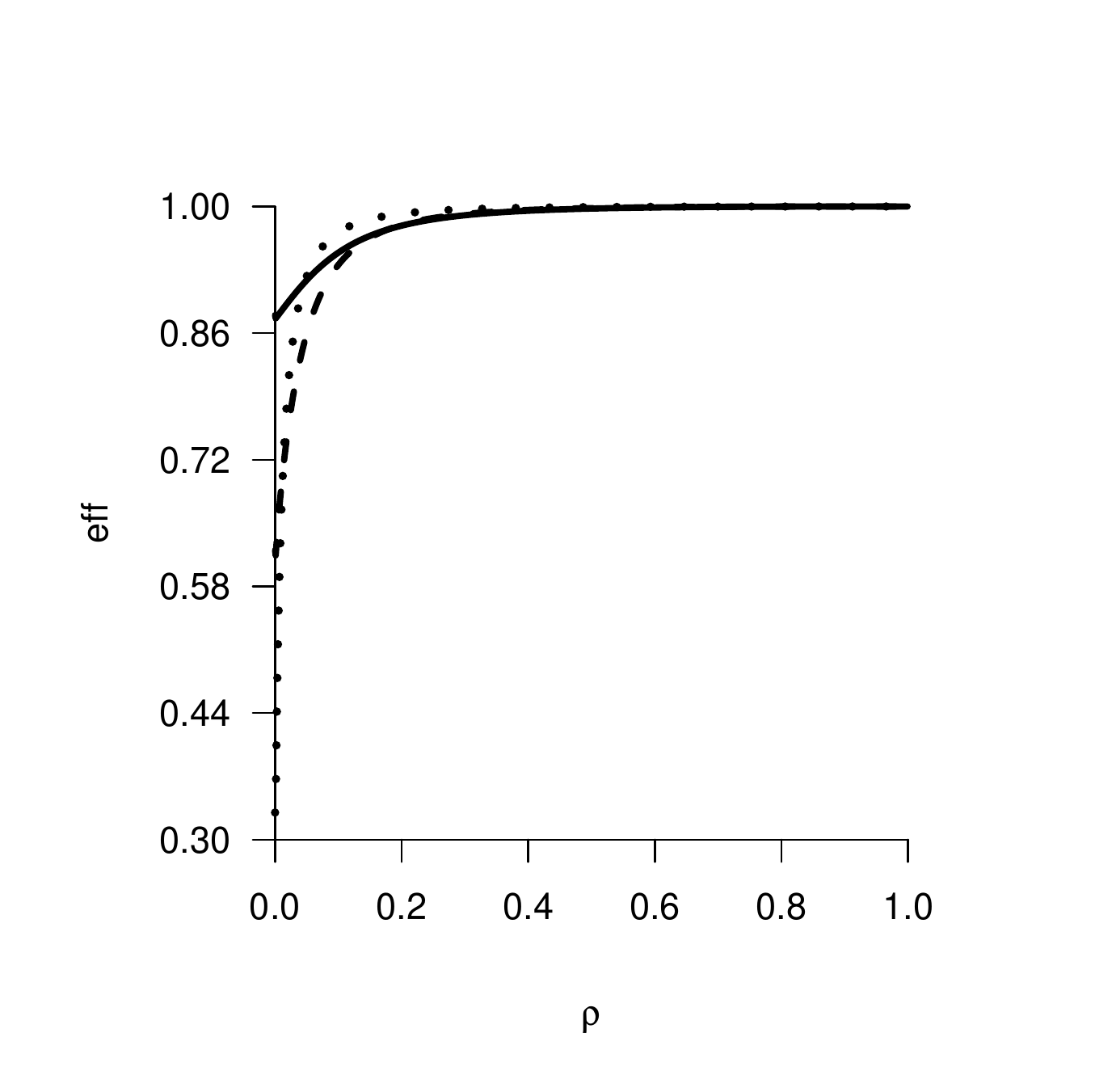}
       \end{minipage}
       \begin{minipage}[]{8.2 cm}
       \centering
       \includegraphics[width=78mm]{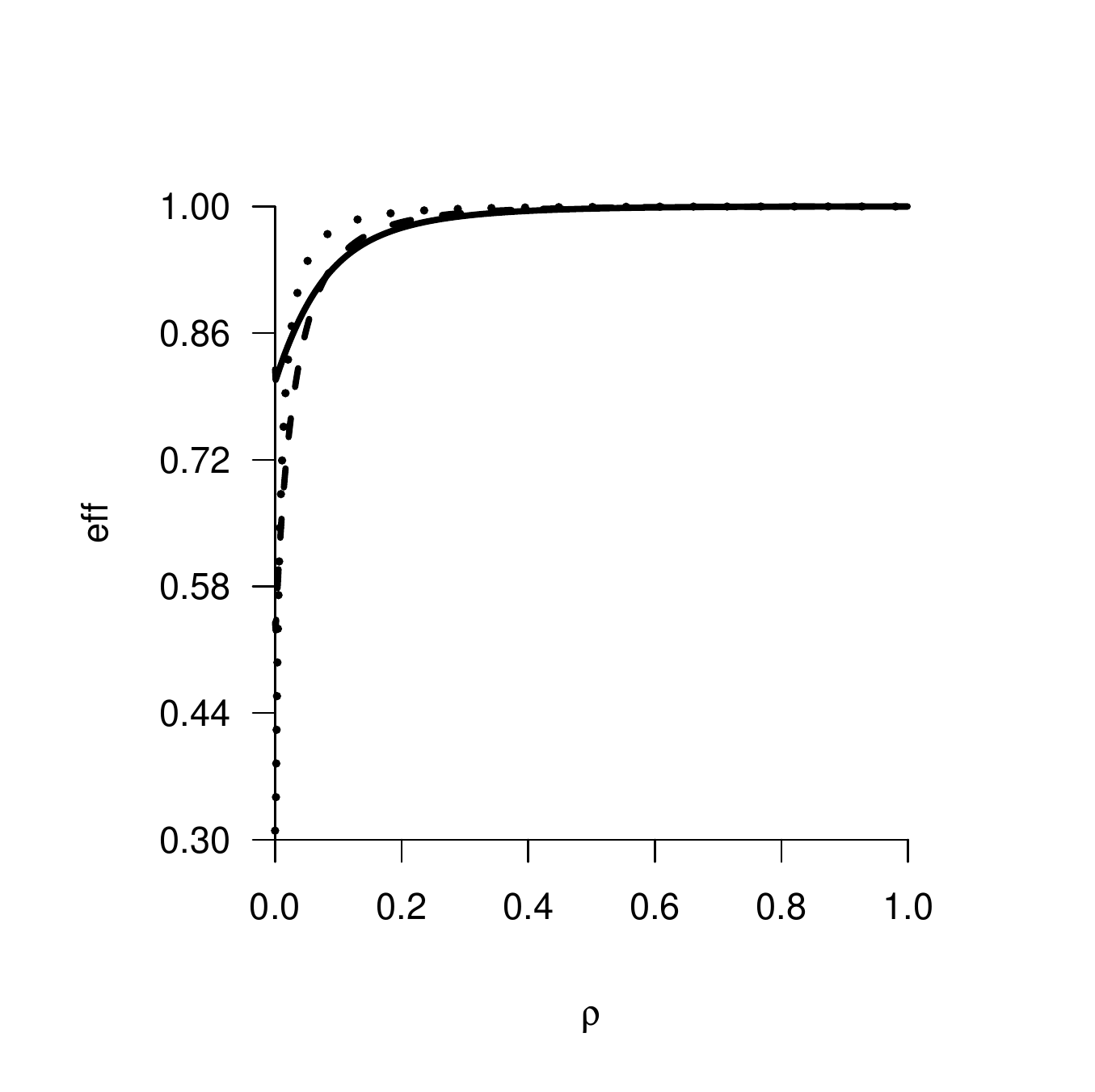}
       \end{minipage}
       \vspace{-1mm}
       \hspace*{5 mm}
       \begin{minipage}[]{6.5 cm}
       Figure 7: Efficiency of minimax-optimal designs for quadratic regression, case~1,
			for $n=10$ (solid line), $n=50$ (dashed line), $n=500$ (dotted line)
       \end{minipage}
       \hspace{15 mm}
       \begin{minipage}[]{6.5 cm}
       Figure 8: Efficiency of minimax-optimal designs for quadratic regression, case~2,
			for $n=10$ (solid line), $n=50$ (dashed line), $n=500$ (dotted line)
       \end{minipage}
    \end{figure}
		The efficiency turns out to be high and increasing with increasing variance parameters for both cases~1 and 2 and all values of the number of individuals ($n=10$, $n=50$, $n=500$).
		
		\section{Discussion}
		
		In this paper we have considered minimax-optimal designs for the IMSE-criterion for the prediction in particular RCR models: linear and quadratic regression. We have assumed the diagonal structure of the covariance matrix of random effects. In this case the IMSE-criterion is increasing with increasing values of all variance parameters. If all variances converge to infinity, the limiting criterion coincides with the IMSE-criterion in fixed effects models and, consequently, the optimal designs in fixed effects models retain their optimality for the prediction. If some of variances are small, the minimax-optimal designs in RCR depend on the number of individuals and differ from the optimal designs in fixed effects models. For some particular cases we have considered the efficiency of the minimax-optimal designs with respect to the locally optimal designs. The efficiency turns out to be high and increase with increasing variance parameters.

\section*{Acknowledgment} This research has been supported by grant SCHW 531/16-1 of the German Research Foundation (DFG).

\bibliographystyle{natbib}
\bibliography{prus4}

\end{document}